\documentclass{article}

\usepackage{amssymb,amscd,amsthm,latexsym}
\usepackage[all]{xy}

\newtheorem{defi}{Definition}[section]
\newtheorem{theo}{Theorem}[section]
\newtheorem{lemma}{Lemma}[section]
\newtheorem{prop}{Proposition}[section]
\newtheorem{coro}{Corollary}[section]
\newtheorem{rem}{Remark}[section]
\newtheorem{ex}{Example}[section]
\begin{document}
\date{}
\author{Dominique Bourn}
\title{The cohomological comparison arising from the associated abelian object}

\maketitle

\begin{abstract}
We make explicit some conditions on a semi-abelian category $\mathbb D$ such that the cohomology group homomorphisms $j^n_A:H^n_{M(\mathbb D/Y)}(A)\rightarrow H^n_{\mathbb D/Y}(A)$, induced by the inclusion $j: Ab\mathbb D\rightarrowtail \mathbb D$ of the abelian objects of $\mathbb D$, are actually isomorphisms. These conditions hold when $\mathbb D$ is the category $Gp$ of groups, and this allows us to give a new insight on the Eilenberg-Mac Lane cohomology of groups. They hold also when $\mathbb D$ is the category $\mathbb D=\mathbb K$-$Lie$ of Lie-algebras.
\end{abstract}

\textbf{Introduction}

\medskip

In a non-abelian context, there are several ways to realize the cohomology groups $H^n_{\mathbb E}(A)$ of a category $\mathbb E$ (provided, let us say, it is finitely complete and exact) with coefficients in an internal abelian group $A$ of $\mathbb E$: by means of simplicial objects in $\mathbb E$, as introduced by Duskin \cite{Du} and Glenn \cite{Gl}, but also by means of internal $n$-groupoids in $\mathbb E$, as in \cite{B1}. When $\mathbb E$ is $Ab$ the category of abelian groups, the equivalence between the category of internal $n$-groupoids in $Ab$ and the category of chain complexes of length $n$ make that these $H^n$ coincide with the Yoneda's $Ext^n$. When $\mathbb E$ is $Gp$ the category of groups, the equivalence between the category of internal $n$-groupoids in $Gp$ and the category of crossed $n$-fold complexes make that these $H^n$ coincide (again see \cite{B1}) with the $Opext^n$ of the Eilenberg-Mac Lane cohomology of groups, via the interpretation (independently) given by Holt \cite{Ho} and Huebschmann \cite{Hu}.

Now, any left exact functor: $U: \mathbb E \rightarrow \mathbb E'$ which preserves the regular epimorphisms provides natural comparison group homomorphisms $U^n_A: H^n_{\mathbb E}(A)\rightarrow H^n_{\mathbb E'}(U(A))$. We shall be especially interested by the following situation. Let $\mathbb D$ be a semi-abelian category \cite{JMT}, as are the categories $Gp$ of groups and $\mathbb K$-$Lie$ of $\mathbb K$-Lie algebras. Now, let $Y$ be an object in $\mathbb D$, and $\mathbb D/Y$ the associated slice category. We determine conditions on $\mathbb D$, such the group homomorphisms $j^n_A:H^n_{M(\mathbb D/Y)}(A)\rightarrow H^n_{\mathbb D/Y}(A)$, induced by the inclusion $j_Y: M(\mathbb D/Y)\rightarrowtail \mathbb D/Y$ of the commutative objects of $\mathbb D/Y$ (induced itself by the inclusion $j: Ab\mathbb D\rightarrowtail \mathbb D$ of the abelian objects of $\mathbb D$) are actually group isomorphisms. These conditions are both global ($\mathbb D$ peri-abelian, see Definition \ref{peri} below) and local (the object $Y$ has  projective dimension $1$). 

These conditions hold, for any object $Y$, in the category $Gp$ of groups and also in the category $\mathbb K$-$Lie$ of $\mathbb K$-Lie algebras. In the case of the category $Gp$, this leads to the rather unexpected observation that, in a way, the categories $Ab$ and $Gp$ are so close to each other (through the ''peri-abelian" connection) that the Eilenberg-Mac Lane cohomology is unable to discriminate them. We thank G.\ Janelidze for our helpful discussions.

The first section of the article is devoted to some recalls about $n$-groupoids. The second one deals with the properties, in the regular Mal'cev context, of the inclusion $j: M\mathbb D\rightarrowtail \mathbb D$ from the commutative objects in $\mathbb D$, and of its extensions to the $n$-groupoids. The third one gives some recalls about the realization of the cohomology groups by means of $n$-groupoids; and the forth one introduces the notion of \emph{peri-abelian} category and provides the final result. 

\section{Internal groupoids and n-groupoids}

We shall suppose all our categories finitely complete. Given the following right hand side commutative square, we denote the kernel equivalence relation
of $f$ by $R[f]$ and the induced map between the kernel equivalences by $R(x)$:
$$\xymatrix@=20pt
{
  R[f] \ar[d]_-{R(x)} \ar@<2ex>[r]^{p_{0}} \ar@<-2ex>[r]_{p_{1}} & X  \ar[d]^-{x} \ar[r]^{f} \ar[l]_{s_0} &
  Y \ar[d]^-{y} \\
  R[f'] \ar@<2ex>[r] \ar@<-2ex>[r] & X' \ar[r]_-{f'} \ar[l] & Y'.
}
$$
An internal groupoid $\underline X_1$ in $\mathbb E$ will be presented (see \cite{B2}) as a reflexive graph $(d_0,d_1):X_1\rightrightarrows X_0$ endowed with an operation $d_2$:
$$\xymatrix@=15pt
{
 R[d_0]^2 \ar@(u,u)[r]^{R(d_2)} \ar@<2ex>[r]^{p_{2}} \ar[r]^{p_{1}} \ar@<-2ex>[r]_{p_{0}} & R[d_0] \ar@(u,u)[r]^{d_2} \ar@<1ex>[r]^{p_{1}} \ar@<-1ex>[r]_{p_{0}} &
X_1  \ar@<2ex>[r]^{d_{1}} \ar@<-1ex>[r]_{d_{0}}  & X_0 \ar[l]_{s_0}
}
$$
making the previous diagram satisfy all the simplicial identities, including the ones concerning the degeneracies. In the set theoretical context, this operation $d_2$ associates the composite $\psi.\phi^{-1}$ with any pair $(\phi,\psi)$ of arrows of $\underline X_1$ with same domain. Any equivalence relation $R$ on an object $X$ in $\mathbb E$ provides an internal groupoid:
$$\xymatrix@=15pt
{
 R[p_0]^2 \ar@(u,u)[r]^{p_3} \ar@<2ex>[r]^{p_{2}} \ar[r]^{p_{1}} \ar@<-2ex>[r]_{p_{0}} & R[p_0] \ar@(u,u)[r]^{p_2} \ar@<1ex>[r]^{p_{1}} \ar@<-1ex>[r]_{p_{0}} &
R  \ar@<2ex>[r]^{p_{1}} \ar@<-1ex>[r]_{p_{0}}  & X \ar[l]_{s_0}
}
$$
which, in some formal circumstances, will be denoted by $\underline R_1$.

Let $Grd\mathbb E$ denote the category of internal groupoids and internal functors in $\mathbb E$, and $()_0:Grd\mathbb E\rightarrow E$ the forgetful functor associating with the groupoid $\underline X_1$ its ''object of objects" $X_0$. This functor is a left exact fibration. Any fibre (above a given object $X$) has a terminal object $\underline{\nabla}_1(X)$ which is the undiscrete relation on the object $X$:
$$
\xymatrix@=20pt
{ 
  X\times X \ar@<2ex>[r]^{p_1} \ar@<-2ex>[r]_{p_0} & X \ar[l]_{s_0}
 }
$$
and an initial object $\underline{\Delta}_1(X)$ which is the discrete equivalence relation on $X$: 
$$
\xymatrix@=20pt
{ 
  X \ar@<2ex>[r]^{1_X} \ar@<-2ex>[r]_{1_X} & X \ar[l]_{1_X}
 }
$$
They produce respectively a right adjoint and a left adjoint to the forgetful functor $()_0$.

An internal functor $\underline f_1: \underline X_1\rightarrow \underline Y_1$ is $()_0$-cartesian if and only if the following square is a pullback in $\mathbb E$, in other words if and only if it is internally fully faithful:
$$\xymatrix@=15pt
{
 X_1 \ar[r]^{f_1} \ar[d]_{(d_0,d_1)} & Y_1 \ar[d]^{(d_0,d_1)}\\
 {X_0\times X_0\;}  \ar[r]_{f_0\times f_0}& Y_0\times Y_0
}
$$
We shall need also the following classical definition:
\begin{defi}\label{fun}
The internal functor $\underline f_1$ is said to be a discrete fibration when any of the following squares is a pullback:
$$\xymatrix@=4pt
{
 X_1 \ar@<1ex>[rrrr]^{d_1} \ar@<-1ex>[rrrr]_{d_0} \ar@<-1ex>[ddd]_{f_1}  & & & & X_0 \ar@<-1ex>[ddd]^{f_0} \ar[llll] \\
 \\
 \\ 
 Y_1  \ar@<1ex>[rrrr]^{d_1} \ar@<-1ex>[rrrr]_{d_0}   && & & Y_0  \ar[llll]
}
$$
\end{defi}

\subsection{Internal n-groupoids} 

They are defined by induction on the integer $n\in \mathbb N$ (\cite{B1}). The category $2$-$Grd\mathbb E$ \linebreak of internal $2$-groupoids is defined as the category of internal groupoids inside the fibres of the fibration $()_0$. We get a forgetful functor $()_{1}:2$-$Grd\mathbb E \rightarrow Grd\mathbb E$ \linebreak associating with any $2$-groupoid $\underline X_2$ its underlying groupoid $\underline X_1$ of $1$-cells. Again, it is a left exact fibration. Now, suppose defined the left exact fibration: $()_{n-2}:(n$-$1)$-$Grd\mathbb E \rightarrow (n$-$2)$-$Grd\mathbb E$, the $n$-groupoids are the internal groupoids inside the fibres of $()_{n-2}$. In this way, an internal $n$-groupoid $\underline X_n$ determines an underlying diagram of reflexive graphs in $\mathbb E$ we shall denote by:
$$
\xymatrix@=15pt{
  \underline X_n : \;\; X_n \ar@<2ex>[r]^-{d_{1}}  \ar@<-2ex>[r]_-{d_0} & X_{n-1}\ar[l]_-{s_0} \ar@<2ex>[r]^-{d_{1}}
  \ar@<-2ex>[r]_-{d_0} & X_{n-2} \ar[l]_-{s_0} \cdots X_1\ar@<2ex>[r]^-{d_{1}}  \ar@<-2ex>[r]_-{d_0} &
  X_0 \ar[l]_-{s_0} }
$$
By $\underline{\nabla}_n$ (respt. $\underline{\Delta}_n$), we denote the right (resp. left) adjoint to the fibration $()_{n-1}$ which gives the terminal (resp. initial) objects in the fibres of $()_{n-1}$. 
To get that, given a $n$-groupoid $\underline X_n$, define by induction the object $X^{\bumpeq}_{n}$ of its "parallel" $n$-cells, from the following string of \emph{kernel equivalence relations}; begin with:
$$
\xymatrix@=20pt
{ 
  X^{\bumpeq}_{1} \ar@<2ex>[r]^{p_1} \ar@<-2ex>[r]_{p_0} & X_{1} \ar[l]_{s_0} \ar[r]^{(d_0,d_1)} &  X_0\times X_0
 }
$$
whose left hand side part determines the upper level of the $2$-groupoid $\underline{\nabla}_{2}\underline X_1$, then construct by induction:
$$
\xymatrix@=20pt
{ 
  X^{\bumpeq}_{n} \ar@<2ex>[r]^{p_1} \ar@<-2ex>[r]_{p_0} & X_{n} \ar[l]_{s_0} \ar[r]^{(d_0,d_1)} &  X^{\bumpeq}_{n-1}
 }
$$
whose left hand side part determines the upper level of the $(n+1)$-groupoid $\underline{\nabla}_{n+1}\underline X_n$. The $(n+1)$-groupoid $\underline{\Delta}_{n+1}\underline X_n$ is given by the discrete equivalence relation on $\underline X_n$.
 
\begin{defi}
We shall say that the $n$-groupoid $\underline X_n$ is pointed when the terminal map $\underline X_n\rightarrow \underline{\nabla}_n(\underline X_{n-1})$ in the fibre of $()_{n-1}$ is split.
\end{defi}

Finally, as an internal groupoid in the fibres of $()_{n-1}$, the $n$-groupoid $\underline X_n$ will be denoted in the following way inside the category $(n$-$1)$-$Grd\mathbb E$:
$$
\xymatrix@=20pt{
  \oint_{n-1} \underline X_n \ar@<2ex>[r]^-{\underline d_{n-1}^1} \ar@<-2ex>[r]_-{\underline d_{n-1}^0} & \underline X_{n-1}, \ar[l] }
$$
where $\oint_{n-1} \underline X_n$ is the integral $(n$-$1)$-groupoid of the $n$-cells of the $n$-groupoid $\underline X_n$:
$$
\xymatrix@=20pt{
  \oint_{n-1} \underline X_n : \;\;  X_n \ar@<2ex>[r]^-{d_1\cdot d_1} \ar@<-2ex>[r]_-{d_0\cdot d_0} &
  X_{n-2} \ar[l]_-{s_0.s_0} \ar@<2ex>[r]^-{d_1} \ar@<-2ex>[r]_-{d_0} & X_{n-3} \ar[l]_-{s_0} \cdots X_1 \ar@<2ex>[r]^-{d_1} \ar@<-2ex>[r]_-{d_0} & X_0. \ar[l]_-{s_0} }
$$
When the category $\mathbb E$ is regular \cite{Ba}, we call \emph{aspherical} a groupoid $\underline X_1$ which is \emph{connected}, namely such that its \emph{core} $(d_0,d_1): X_1\rightarrow X_0\times X_0$ is a regular epimorphism, and which has its object of objects $X_0$ with global support, namely such that terminal map $X_0\rightarrow 1$ is a regular epimorphism. A $n$-groupoid $\underline X_n$ is \emph{aspherical} when it is \emph{connected}, namely such that its \emph{core}  $(d_0,d_1): X_n\rightarrow X^{\bumpeq}_{n-1}$ is a regular epimorphism, and which has its underlying $(n$-$1)$-groupoid $\underline X_{n-1}$ aspherical.

\begin{defi}
A $n$-functor $\underline f_{n}: \underline X_{n} \rightarrow \underline Y_{n}$ is said to be a discrete fibration when any of the following squares (at the highest level) is a pullback:
$$\xymatrix@=4pt
{
 X_n \ar@<1ex>[rrrr]^{d_1} \ar@<-1ex>[rrrr]_{d_0} \ar@<-1ex>[ddd]_{f_n}  & & & & X_{n-1} \ar@<-1ex>[ddd]^{f_{n-1}} \ar[llll] \\
 \\
 \\ 
 Y_n  \ar@<1ex>[rrrr]^{d_1} \ar@<-1ex>[rrrr]_{d_0}   && & & Y_{n-1}  \ar[llll]
}
$$
\end{defi}

\begin{prop}\label{df}
Suppose $n\geq 1$ and the $(n+1)$-functor $\underline f_{n+1}: \underline X_{n+1} \rightarrow \underline Y_{n+1}$ \linebreak such that the underlying $n$-functor $\underline f_n$ is $()_{n-1}$-cartesian.  Then  $\underline f_{n+1}$ is a discrete fibration if and only if it is $()_{n}$-cartesian. Suppose moreover the category $\mathbb E$ regular, the $(n+1)$-groupoid $\underline X_{n+1}$ asherical and the $(n+1)$-groupoid $\underline Y_{n+1}$ connected, then there is a kind of converse: when $\underline f_{n+1}$ is at the same time $()_{n}$-cartesian and a discrete fibration, then the underlying $n$-functor $\underline f_{n}$ is $()_{n-1}$-cartesian.
\end{prop}
\proof
Consider the following diagram:
$$
\xymatrix@=15pt
{ X_{n+1} \ar[rr]^{f_{n+1}}\ar[d]^{(d_0,d_1)} \ar@(l,l)[dd]_{d_0} && Y_{n+1} \ar[d]_{(d_0,d_1)} \ar@(r,r)[dd]^{d_0} \\
  X^{\bumpeq}_{n} \ar@<-1ex>[d]_{p_0} \ar@<1ex>@{..>}[d]^{p_1} \ar[rr]^{f^{\bumpeq}_{n}}& & Y^{\bumpeq}_{n} \ar@<-1ex>[d]_{p_0} \ar@<1ex>@{..>}[d]^{p_1}
   \\ 
  X_n \ar[rr]^{f_{n}} \ar[d]_{(d_0,d_1)} && Y_n \ar[d]^{(d_0,d_1)} \\
  X^{\bumpeq}_{n-1} \ar[rr]_{f^{\bumpeq}_{n-1}} && Y^{\bumpeq}_{n-1}
 }
$$
The fact that the $n$-functor $\underline f_n$ is $()_{n-1}$-cartesian means that the lower square is a pullback. Then the middle square is a pullback. Accordingly the bow quadrangle is a pullback (which precisely means that $\underline f_{n+1}$ is a discrete fibration) if and only if the the upper square is a pullback (which means that $\underline f_{n+1}$ is $()_{n}$-cartesian).

Suppose now that $\mathbb E$ is regular and that the bow quadrangle and the upper square are pullbacks. If moreover the two upper vertical arrows are regular epimorphisms (which means that the $(n+1)$-groupoids $\underline X_{n+1}$ and $\underline Y_{n+1}$ are connected), then the middle square is a pullback. When moreover the left hand side lower vertical arrow is a regular epimorphism (which follows from the fact that $\underline X_{n+1}$ is aspherical), then the lower square is a pullback (by the Barr-Kock theorem), and the underlying $n$-functor $\underline f_{n}$ is $()_{n-1}$-cartesian. 
\endproof

\subsection{The regular Mal'cev context}

A category $\mathbb D$ is a Mal'cev category when any reflexive relation is an equivalence relation, see \cite{CLP} and \cite{CPP}. Equivalently  $\mathbb D$ is a Mal'cev category when any reflexive sub-graph of an internal groupoid is a groupoid. In this context, any internal category is a groupoid. The categories $Gp$ of groups, $Rg$ of non unitary commutative rings and $\mathbb K$-$Lie$ of $\mathbb K$-Lie algebras are Mal'cev categories. When $\mathbb D$ is a Mal'cev category, then any slice category $\mathbb D/Y$ is still a Mal'cev category. 

It appears that the context of Mal'cev categories particularly fits with the notion of commutator of equivalence relations \cite{BG}, \cite{Po}. It is then possible to define an object $X$ in $\mathbb D$ as being \emph{commutative} when we have $[\nabla X,\nabla X]=0$, i.e. when the commutator $[\nabla X,\nabla X]$ is trivial, or, equivalently, when the object $X$ is equipped with a (unique possible) Mal'cev operation $p:X\times X \times X \rightarrow X$. We shall denote by $M\mathbb D$ the subcategory of the \emph{commutative} objects $X$ in $\mathbb D$ and by $j: M\mathbb D \rightarrowtail \mathbb D$ the inclusion functor. A morphism $f: X\rightarrow Y$ is central when we have $[R[f],\nabla X]=0$. It has an abelian kernel relation (or is commutative as an object in the slice category $\mathbb D/Y$) when we have $[R[f],R[f]]=0$. 

In $Gp$, a commutative object is an abelian group. Given any group $Y$, the slice category $Gp/Y$ is still a Mal'cev category, and an object of $Gp/Y$, namely a group homomorphism $X\rightarrow Y$, is commutative in $Gp/Y$ if and only if its kernel is abelian. Given any finitely complete category $\mathbb E$, the fibres of the fibration $()_0:Grd\mathbb E\rightarrow E$ are Mal'cev categories, so that there is a natural notion of \emph{commutative} (actually we say \emph{abelian}) groupoid. Recall that, when $\mathbb D$ is a Mal'cev category, then any internal groupoid in $\mathbb D$ is \emph{abelian}.

Recall a category $\mathbb D$ is regular \cite{Ba} when the regular epimorphisms are stable under pullback and any effective equivalence relation admits a quotient. It is exact when, moreover, any equivalence relation is effective. Any variety of Universal Algebra is exact, and in particular the category $Gp$. We shall need the intermediate notion of \emph{efficiently regular} category \cite{B10}: 

\begin{defi}
\emph{A regular category $\mathbb D$ is said to be \emph{efficiently regular} when any equivalence relation $T$ on an
object $X$ which is a subobject $j : T\rightarrowtail R[f]$ of an effective equivalence relation $R[f]$ on $X$ by an effective
monomorphism in $\mathbb C$ (which means that $j$ is the equalizer of some pair of maps in $\mathbb C$) is itself
effective.}
\end{defi}

Any exact category is efficiently regular. The category $GpTop$ (resp. $AbTop$) of topological (resp. topological abelian) groups is efficiently regular, but not exact. More generally any category $Top^{\mathbb T}$ of topological protomodular algebras (where $\mathbb T$ is a protomodular theory) is efficiently regular: it is a regular category according to \cite{BC}, and clearly an equivalence relation $T$ on $X$ is effective if and only if the object $T$ is endowed with the topology induced by the topological product, which is the case when $j: T\rightarrowtail R[f]$ is an effective monomorphism. When $\mathbb E$ is efficiently regular, such is any slice category $\mathbb E/Y$, and any fibre of the fibration $()_0:Grd\mathbb E \rightarrow \mathbb E$. The main fact in an efficiently regular Mal'cev category $\mathbb D$ is that any commutative object $X$ has a \emph{direction} which is given by the following diagram \cite{B10}:
$$
\xymatrix@C=20pt@R=25pt{
  X\times X\times X \ar@<1ex>[r]^-{p_2} \ar@<-1ex>[r]_-{(p_0,p)} \ar@<-.5ex>[d]_-{p_0} &
  X\times X  \ar@<-.5ex>[d]_-{p_0} \ar@{-->>}[r]^-{\nu_X} & A \ar@<-.5ex>[d] \\
  X\times X \ar@<1ex>[r]^-{p_{1}}  \ar@<-1ex>[r]_-{p_0} \ar@<-.5ex>[u] & X \ar[r] \ar@<-.5ex>[u] & 1 \ar@<-.5ex>@{.>}[u]. }
$$
where the map $\nu_X$ is the quotient of the upper horizontal equivalence relation (where the maps $p_i$ are the product projections, and the map $p$ is the Mal'cev operation which makes $X$ commutative). Any commutative square in this diagram is a pullback. When moreover $X$ has a global support, the section $s_0$ of $p_0$ can be extented to he quotient (by the dotted arrow), $A$ becomes an internal abelian group in $\mathbb D$ and the upward right hand side square becomes a pushout. For the same general reasons, any aspherical groupoid $\underline X_1$, being abelian in the Mal'cev context, admits a \emph{direction}, namely  an internal abelian group $A$ in $\mathbb D$ such that the following downward square is a pullback, while, at the same time, the upward square is a pushout:
$$\xymatrix@=20pt
{
 X^{\bumpeq}_{1} \ar@{->>}[r]^{} \ar@<1ex>[d]^{p_{1}} \ar@<-1ex>[d]_{p_{0}} & A \ar@<1ex>[d]^{} \\ 
 X_{1}   \ar@{->>}[r]_{} \ar[u] \ar@{->>}[d]_{(d_0,d_1)} & 1 \ar[u]\\
 X_0\times X_0
}
$$
More generally, any aspherical $n$-groupoid $\underline X_n$ admits a \emph{direction} $A$ given by the same kind of diagram (see \cite{BR}):
$$\xymatrix@=20pt
{
 X^{\bumpeq}_{n} \ar@{->>}[r]^{} \ar@<1ex>[d]^{p_{1}} \ar@<-1ex>[d]_{p_{0}} & A \ar@<1ex>[d]^{} \\ 
 X_{n}   \ar@{->>}[r]_{} \ar[u] \ar@{->>}[d]_{(d_0,d_1)} & 1 \ar[u]\\
 X^{\bumpeq}_{n-1}
}
$$
Any $n$-functor $\underline f_n: \underline X_n\rightarrow \underline Y_n$ between two aspherical $n$-groupoids produces a group homomorphism between their respective directions.  

\begin{rem}\label{rema}
This group homomorphism is an isomorphism if and only if the $n$-functor $\underline f_n$ is $()_{n-1}$-cartesian, see \cite{B1}.
\end{rem}

Finally when $\mathbb D$ is moreover finitely cocomplete, the inclusion functor $j:M\mathbb D \rightarrowtail \mathbb D$ admits a left adjoint (see \cite{B12}) which will be denoted by $M$. 

\section{The reg-epi reflections}

Now, more generally, we shall consider a regular Mal'cev category $\mathbb D$ and $j:\mathbb C \rightarrowtail \mathbb D$ a full replete inclusion  which admits a reg-epi reflection $I:\mathbb D \rightarrow \mathbb C$. Recall:

\begin{defi}
Let $j:\mathbb C \rightarrowtail \mathbb D$ be a full replete inclusion and $\mathbb D$ a regular category. We shall say that the reflection $I:\mathbb D \rightarrow \mathbb C$ is a reg-epi reflection when any projection $\eta_X:X\twoheadrightarrow IX$ is a regular epimorphism.
\end{defi}

In the regular context, being a reg-epi reflection is equivalent to saying that $\mathbb C$ is stable under subobjects. So, when $\mathbb D$ is a finitely cocomplete regular Mal'cev category, the reflection  $M:\mathbb D \rightarrow M\mathbb D$ to the subcategory $M\mathbb D$ of the commutative objects in $\mathbb D$ is a reg-epi reflection. On the other hand, when $I$ is a reg-epi reflection, a map $h:C\rightarrow C'$ in $\mathbb C$ is a regular epimorphism in $\mathbb C$ if and only if it is a regular epimorphism in $\mathbb D$, so that the inclusion $j$ preserves regular epimorphisms.

Any result of the last part of this section will be borrowed from \cite{BR2}. The reg-epi reflections have strong left exact properties, among which the preservation of internal groupoids: 
\begin{prop}\label{pb}
When $\mathbb D$ is a regular Mal'cev category, any reg-epi reflection $I$ preserves the pullbacks of split epimorphims along split epimorphisms. Accordingly it preserves the kernel equivalence relations of split epimorphisms, and the image $I(\underline X_1)$ of any internal groupoid $\underline X_1$ is an internal groupoid.
\end{prop}

However the kernel equivalence relation of any map is not preserved by $I$ in general, and we shall need to identify a certain class of maps having this property. Following \cite{JK}, we have the following:
\begin{defi}
A map $f:X\rightarrow Y$ in $\mathbb D$ is said to be $I$-trivial when the following square is a pullback:
$$
\xymatrix@=20pt
{ 
  X \ar@{->>}[r]^{\eta_X} \ar[d]_{f} & IX \ar[d]^{If}\\
  Y \ar@{->>}[r]_{\eta_Y} & IY
 }
$$ 
\end{defi}

Any map in $\mathbb C$ is $I$-trivial. The isomorphisms are $I$-trivial, the $I$-trivial maps are stable under composition and such that, when $g.f$ and $g$ are $I$-trivial, then $f$ is $I$-trivial. Also $I$-trivial maps are stable under those pullbacks which are preserved by the reflection $I$. On the other hand, a $I$-trivial map $f$ is certainly $I$-cartesian, namely universal among the maps above $I(f)$. 

\begin{prop}\label{n4}
Suppose $\mathbb D$ is a regular Mal'cev category, $I$ is a reg-epi reflection and $f$ is $I$-trivial. Then we have $I(R[f])\simeq R[I(f)]$, and the maps $p_i:R[f]\rightarrow X$ and $s_0:X\rightarrowtail R[f]$ are still $I$-trivial. In other words, the reflection $I$ preserves the kernel equivalence relations of the $I$-trivial maps.
\end{prop}

\begin{prop}\label{n2}
Suppose $\mathbb D$ is an efficiently regular Mal'cev category and $I$ a reg-epi reflection. Then a regular epimorphism is $I$-trivial if and only if it is the pullback of some map of the subcategory $\mathbb C$. The $I$-trivial regular epimorphisms coincide with the $I$-cartesian regular epimorphisms, they are stable under pullbacks along any map, and these pullbacks are preserved by $I$.
\end{prop}

Again following \cite{JK}, we shall need also the following extension of the class of $I$-trivial maps:

\begin{defi}
Suppose $\mathbb D$ is a regular Mal'cev category and $I$ a reg-epi reflection . We call a map $f:X\rightarrow Y$ $I$-normal when the projection $p_0:R[f]\rightarrow X$ is $I$-trivial.
\end{defi}

Then, by Proposition \ref{n4}, any $I$-trivial map is $I$-normal. Considering any finitely cocomplete efficiently regular Mal'cev category $\mathbb D$ and the reg-epi reflection  $M: \mathbb D\rightarrow M\mathbb D$, a regular epimorphism $f:X\twoheadrightarrow Y$ is $M$-normal if and only if it is central is the classical sense (i.e. $[R[f],\nabla X]=0$) \cite{JK}, \cite{BR2}. Finally we get:

\begin{prop}\label{n3}
Suppose $\mathbb D$ is a regular Mal'cev category and $I$ a reg-epi reflection. If the regular epimorphism $f:X\twoheadrightarrow Y$ is $I$-normal, then it is $I$-trivial if and only if its kernel equivalence relation is preserved by $I$ (i.e. $R[I(f)]\simeq I(R[f])$). Consequently a split epimorphism $f$ is $I$-trivial if and only if it is $I$-normal.
\end{prop}

We shall need also the following extension of the last assertion:
\begin{prop}\label{spec}
Suppose $\mathbb D$ is a regular Mal'cev category and $I$ a reg-epi reflection. Any split epimorphism $t$ between two $I$-normal maps $f$ and $f'$ in the slice category $\mathbb D/Y$:
$$\xymatrix@=7pt
{ 
X  \ar[rr]_{t} \ar[ddr]_f  && X'\ar@<-1ex>[ll]_{s} \ar[ddl]^{f'}\\
\\ 
& Y 
}
$$
is $I$-trivial. 
\end{prop}

\subsection{The reg-epi reflections $I: n$-$Grd\mathbb D \rightarrow n$-$Grd\mathbb C$}\label{sec2}

Suppose $\mathbb D$ is a regular Mal'cev category and $I: \mathbb D\rightarrow \mathbb C$ a reg-epi reflection. We noticed that $I$ preserves the internal groupoids. Consequently, for any $n\in\mathbb N$, the functor $I$ extends naturally to a functor: $n$-$Grd\mathbb D \rightarrow n$-$Grd\mathbb C$ (for sake of simplicity still denoted by $I$) which, again, is a reg-epi reflection. We are now going to investigate it more precisely.\\ 

From now on, we shall suppose $\mathbb D$ is an efficiently regular Mal'cev category (in order to be able to define the direction of aspherical groupoids) and $I:\mathbb D \rightarrow \mathbb C$ is a reg-epi reflection. We shall assume moreover that the reflection $I:\mathbb D \rightarrow \mathbb C$ :\\
1) induces an isomorphism $Ab\mathbb C \simeq Ab\mathbb D$, where $Ab\mathbb D$ denotes the category of abelian groups in $\mathbb D$, or equivalently the category of pointed commutative objects in $\mathbb D$\\
2) preserves the groupoids $\underline{\nabla}_1(X)$ (i.e. is such that $I(X\times X)\simeq IX\times IX$).\\  

The main example we have in mind satisfying those conditions is the reg-epi reflection $M: \mathbb D \rightarrow M\mathbb D$ towards the commutative objects in a finitely cocomplete efficiently regular Mal'cev category $\mathbb D$. It obviously satisfies Condition 1; it is shown in \cite{BR2} that, under the assumption that $\mathbb D$ is efficiently regular, it satisfies also Condition 2.

\medskip 

What is important for us, with the condition $Ab\mathbb C\simeq Ab\mathbb D$, is that any aspherical groupoid $\underline X_1$ has its core $(d_0,d_1): X_1\twoheadrightarrow X_0\times X_0$ $I$-normal:
\begin{lemma}
Under the assumptions precised at the beginning of this section, any aspherical groupoid $\underline X_1$ has its core $X_1\twoheadrightarrow X_0\times X_0$ $I$-normal.
\end{lemma}
\proof
Indeed, let $A$ be its direction. As an abelian group in $\mathbb D$, it lies also in $\mathbb C$. According to the Proposition \ref{n2} the following upper pullbacks show that $p_0$ is $I$-trivial and, consequently, the core $(d_0,d_1): X_1\twoheadrightarrow X_0\times X_0$ is $I$-normal:
$$\xymatrix@=20pt
{
 X^{\bumpeq}_{1} \ar@{->>}[r]^{} \ar@<1ex>[d]^{p_{1}} \ar@<-1ex>[d]_{p_{0}} & A \ar@<1ex>[d]^{} \\ 
 X_{1}   \ar[r]_{} \ar[u] \ar@{->>}[d]_{(d_0,d_1)} & 1 \ar[u]\\
 X_0\times X_0
}
$$
\endproof

\noindent The functor $I$ preserves the aspherical groupoids, since $I$ preserves the regular epimorphism and is such that we have $I(X\times X)=I(X)\times I(X)$. Here are our first problematical observations:

There is no reason why, in general, \emph{the functor $I$ would preserve the direction of the aspherical groupoids} (which would be equivalent to saying that the functor $\underline X_1\twoheadrightarrow I\underline X_1$ is $()_0$-cartesian), since the objects $X_1^{\bumpeq}$ are not preserved by $I$. And no reason why the image of any aspherical $n$-groupoid, $2\leq n$,  would be aspherical. However we can get the following precisions:
\begin{prop}
Let $\underline X_{2}$ be an aspherical $2$-groupoid, then $d_0:X_2\rightarrow X_1$ is $I$-trivial and $\eta_2:\underline X_{2}\rightarrow I\underline X_{2}$ is a discrete fibration. 
\end{prop}
\proof
Since the $1$-groupoids $\oint_1 \underline X_2$ and $\underline X_1$ are both aspherical, the cores \linebreak $(d_0.d_0,d_1.d_1): X_2\rightarrow X_0\times X_0$ and $(d_0,d_1): X_1\rightarrow X_0\times X_0$ are $I$-normal. Thus the split epimorphism $d_0:X_2\rightarrow X_1$, which makes commute these cores, is $I$-trivial by Proposition \ref{spec}, and $\eta_2:\underline X_{2}\rightarrow I\underline X_{2}$ is a discrete fibration.
\endproof

\begin{prop}
Let $\underline X_{2}$ be an aspherical $2$-groupoid, then its core $(d_0,d_1):X_2\rightarrow X^{\bumpeq}_{1}$ is $I$-trivial. 
\end{prop}
\proof
We have $d_0=p_0.(d_0,d_1): X_2\rightarrow X^{\bumpeq}_{1} \rightarrow X_1$. The map $d_0:X_2\rightarrow X_1$ is $I$-trivial. This is equally the case for $p_0$
$$
\xymatrix@=30pt
{ 
  X^{\bumpeq}_{1} \ar@<2ex>[r]^{p_1} \ar@<-2ex>[r]_{p_0} & X_{1} \ar[l]_{s_0} \ar[r]^{(d_0,d_1)} &  X_0\times X_0
 }
$$
since $(d_0,d_1):X_1\rightarrow X_0\times X_0$ is $I$-normal. Accordingly $(d_0,d_1):X_2\rightarrow X^{\bumpeq}_{1}$ is $I$-trivial.
\endproof

\begin{prop}\label{dfn}
Let $\underline X_{n}$ be an aspherical $n$-groupoid, $n\geq 2$, then the $n$-functor $\eta_n:\underline X_{n}\rightarrow I\underline X_{n}$ is a discrete fibration and its core $(d_0,d_1):X_n\rightarrow X^{\bumpeq}_{n-1}$ is $I$-trivial. 
\end{prop}
\proof
By induction. We saw it is true when $n=2$. Suppose it is true as far as the level $n$-$1$. The $(n$-$1)$-groupoids $\oint_{n-1} \underline X_n$  and $\underline X_{n-1}$ being aspherical, the maps $(d_0,d_1):X_{n}\rightarrow X^{\bumpeq}_{n-2}$ and $(d_0,d_1):X_{n-1}\rightarrow X^{\bumpeq}_{n-2}$ are $I$-trivial by the inductive assumption. Accordingly, the map $d_0:X_n\rightarrow X_{n-1}$ is $I$-trivial and $\eta_n:\underline X_{n}\rightarrow I\underline X_{n}$ is a discrete fibration. Moreover we have $d_0=p_0.(d_0,d_1): X_n\rightarrow X^{\bumpeq}_{n-1} \rightarrow X_{n-1}$. The map $d_0:X_n\rightarrow X_{n-1}$ is $I$-trivial. This is equally the case for the map $p_0$: 
$$
\xymatrix@=30pt
{ 
  X^{\bumpeq}_{n-1} \ar@<2ex>[r]^{p_1} \ar@<-2ex>[r]_{p_0} & X_{n-1} \ar[l]_{s_0} \ar[r]^{(d_0,d_1)} &  X^{\bumpeq}_{n-2}
 }
$$ 
since $(d_0,d_1):X_{n-1}\rightarrow X^{\bumpeq}_{n-2}$ is $I$-trivial by the inductive assumption. Accordingly $(d_0,d_1):X_n\rightarrow X^{\bumpeq}_{n-1}$ is $I$-trivial.
\endproof 

So the only level where the core of an aspherical $n$-groupoid is not necessarily $I$-trivial is the level 1.

\begin{defi}
A groupoid $\underline X_1$ will be said to be $I$-specific when it is aspherical and such that its core $(d_0,d_1):X_1\rightarrow X_0\times X_0$ is $I$-trivial. More generally, a $n$-groupoid $\underline X_n$ will be said to be $I$-specific when it is aspherical and such that $\underline X_{n-1}$ is $I$-specific.
\end{defi}

We are now going to show that this special class of aspherical $n$-groupoids has a "better" behaviour with respect to $I$:

\begin{ex}\label{pointed}
Since any aspherical groupoid $\underline X_1$ has its core $X_1\twoheadrightarrow X_0\times X_0$ $I$-normal, then any pointed groupoid is $I$-specific, following Proposition \ref{n3}. Actually this last point is equivalent to the condition $Ab\mathbb C\simeq Ab\mathbb D$.
\end{ex}
In presence of the assumption 2) about the functor $I$, the groupoid $\underline X_1$ is $I$-specific if and only if the functor $\eta_1:\underline X_{1}\rightarrow I\underline X_{1}$ is $()_{0}$-cartesian or equivalently, (following Proposition \ref{n3}), if and only if we have $\underline{\nabla}_2(I(\underline X_1))=I(\underline{\nabla}_2(\underline X_1))$.

\begin{prop}\label{undis}
Let $\underline X_{n}$ be an aspherical $n$-groupoid, then it is $I$-specific if and only if, for all $1\leq k \leq n$, we have $\underline{\nabla}_{k+1}(I(\underline X_k))=I(\underline{\nabla}_{k+1}(\underline X_k))$.
\end{prop}
\proof
By induction. The equivalence holds at level 1. Suppose it holds as far as the level $n$-$1$. Then consider the following diagram where all the left hand side horizontal arrows are the maps $\eta$:
$$\xymatrix@=20pt
{
 X^{\bumpeq}_{n} \ar@{->>}[r]^{} \ar@<1ex>[d]^{p_{1}} \ar@<-1ex>[d]_{p_{0}} \ar@<2ex>[rr]^{(\eta_{X_n})^{\bumpeq}} &I(X^{\bumpeq}_{n}) \ar[r]_{\epsilon_{n}} \ar@<1ex>[d]^{I(p_{1})} \ar@<-1ex>[d]_{I(p_{0})} & I(X_n)^{\bumpeq} \ar@<1ex>[d]^{p_{1}} \ar@<-1ex>[d]_{p_{0}} \\ 
 X_{n}   \ar@{->>}[r]^{\eta_{X_n}} \ar[u] \ar[d]_{(d_0,d_1)} & I(X_{n})   \ar[r]_{1_{I(X_n)}} \ar[u] \ar[d]_{I((d_0,d_1))} & I(X_{n}) \ar[u] \ar[d]_{(d_0,d_1)}\\
 X^{\bumpeq}_{n-1}  \ar@{->>}[r]_{}\ar@<-2ex>[rr]_{(\eta_{X_{n-1}})^{\bumpeq}}  & I(X^{\bumpeq}_{n-1})\ar[r]^{\epsilon_{n-1}} & I(X_{n-1})^{\bumpeq}
}
$$
Suppose $\underline X_{n}$ $I$-specific. Then $\underline X_{n-1}$ is $I$-specific, and  the inductive assumption ($\underline{\nabla}_{k+1}(I(\underline X_k))=I(\underline{\nabla}_{k+1}(\underline X_k))$, $1\leq k \leq n-1$) implies that $\epsilon_{n-1}$ is an isomorphism. Since the middle vertical diagram is a kernel equivalence relation ($(d_0,d_1):X_{n}\rightarrow X^{\bumpeq}_{n-1}$ being $I$-trivial by Proposition \ref{dfn}), the factorization $\epsilon_{n}$ is an isomorphism, and we get $\underline{\nabla}_{n+1}(I(\underline X_n))=I(\underline{\nabla}_{n+1}(\underline X_n))$ 

Conversely, suppose we have $\underline{\nabla}_{k+1}I(\underline X_k))=I(\underline{\nabla}_{k+1}(\underline X_k))$ for all $1\leq k \leq n$. The inductive assumption is satisfied as far as level $n$-$1$. Accordingly $\underline X_{n-1}$, and thus $\underline X_{n}$, is $I$-specific
\endproof

\begin{coro}
Given any $I$-specific $n$-groupoid $\underline X_n$, the groupoid $I(\underline X_n)$ is aspherical.
\end{coro}

\begin{prop}\label{ff}
Let $\underline X_{n}$ be an aspherical $n$-groupoid, then it is $I$-specific if and only if $\eta_n:\underline X_{n}\rightarrow I\underline X_{n}$ is $()_{n-1}$-cartesian. Actually it is $I$-specific if and only if $\eta_1:\underline X_{1}\rightarrow I\underline X_{1}$ is $()_{0}$-cartesian.
\end{prop}
\proof
Recall that when $\underline X_{n}$ is aspherical, the map $\eta_n:\underline X_{n}\rightarrow I\underline X_{n}$ is a discrete fibration. The direct proof will be given by induction. This is true for $n=1$, since $(d_0,d_1):X_1\rightarrow X_0\times X_0$ $I$-trivial means precisely $\eta_1:\underline X_{1}\rightarrow I\underline X_{1}$ $()_0$-cartesian. Suppose the result holds as far as the level $n$-$1$. If $\underline X_{n}$ is $I$-specific, then $\underline X_{(n-1)}$ is $I$-specific. By the inductive assomption $\eta_{n-1}:\underline X_{n-1}\rightarrow I\underline X_{n-1}$ is $()_{n-2}$-cartesian, and by Proposition \ref{df} $\eta_n:\underline X_{n}\rightarrow I\underline X_{n}$ is $()_{n-1}$-cartesian. Again the converse is true for $n=1$. Suppose it holds as far as the level $n$-$1$. \linebreak If $\eta_n:\underline X_{n}\rightarrow I\underline X_{n}$ is $()_{n-1}$-cartesian, then, since $\eta_n:\underline X_{n}\rightarrow I\underline X_{n}$ is a discrete fibration and the $n$-groupoids $\underline X_{n}$ and $I(\underline X_{n})$ are aspherical, the $(n-1)$-functor $\eta_{n-1}:\underline X_{n-1}\rightarrow I\underline X_{n-1}$ is $()_{n-2}$-cartesian by Proposition \ref{df}; and by the inductive assumption the $(n$-$1)$-groupoid $\underline X_{n-1}$, and thus the $n$-groupoid $\underline X_{n}$, are $I$-specific. The last point comes from the fact that any $\eta_k:\underline X_{k}\rightarrow I\underline X_{k}$ is a discrete fibration by Proposition \ref{df}.
\endproof

As a consequence, by Remark \ref{rema}, we then get what we were aiming to:

\begin{theo}
Given any $I$-specific $n$-groupoid $\underline X_n$, the groupoids $I(\underline X_n)$ and $\underline X_n$ have same direction.
\end{theo}

\section{The cohomology groups with coefficients in $A$}

Let $\mathbb E$ be an efficiently regular category, and $A$ be an internal abelian group in $\mathbb E$. Then, as usual, define $H^0_{\mathbb E}(A)$ as the abelian group $Hom_{\mathbb E}(1,A)$ of global sections of $A$, and $H^1_{\mathbb E}(A)$ as the abelian group of $A$-torsors, namely objects $X$ with global support endowed with a simply transitive action of the group $A$. There are several equivalent ways to realize the groups $H^n_{\mathbb E}(A)$, $n\geq 2$: by means of simplicial objects, as introduced by Duskin \cite{Du} and Glenn \cite{Gl}, but also by means of $n$-groupoids, as in \cite{B1}. In this last case, when $\mathbb E=\mathbb A/C$, where $\mathbb A$ is an abelian category, these $H^n_{\mathbb A/C}(A)$ coincide with the Yoneda $Ext^n(C,A)$ via the equivalence between the extensions of length $n$ in $\mathbb A$ and the aspherical $n$-groupoids in $\mathbb A$. When $\mathbb E=Gp/Y$, these $H^n_{Gp/Y}(A)$ coincide with the $Opext^n(Y,A)$ of the Eilenberg-Mac Lane cohomology of groups, via the interpretation given by Holt \cite{Ho} and Huebschmann \cite{Hu}, and the equivalence between the crossed $n$-fold extensions and the aspherical $n$-groupoids in $Gp$, see \cite{BH} and \cite{B1}. When $\mathbb E=\mathbb D$ is supposed to be a Mal'cev category, the abelian group $H^{n+1}_{\mathbb D}(A)$ appears to be made of the component classes of the aspherical $n$-groupoids with direction $A$, see \cite{BR}, while the group $H^1_{\mathbb D}(A)$ is made of the component classes of the commutative objects with direction $A$. The addition is given by a tensor product analogous to the tensor product of torsors. The inverse of the class of $\underline X_n$ in $H^{n+1}_{\mathbb D}(A)$ is given by the class of its dual $\underline X_n^{op}$ which is defined by inverting the role of the domain and the codomain maps at the last level, the one of the $n$-cells.

\subsection{The cohomological comparison}

Now, suppose we are under the assumptions of section \ref{sec2}, and let $A$ be in  $Ab\mathbb C = Ab\mathbb D$. Then clearly the groups $H^0_{\mathbb C}(A)$ and $H^0_{\mathbb D}(A)$, whose elements are the sections (in $\mathbb C$ and $\mathbb D$ respectively) of the terminal map $A\rightarrow 1$, are the same. Moreover, since the inclusion functor $j:\mathbb C\rightarrowtail \mathbb D$ is left exact and, having a reg-epi reflection $I$, preserves the regular epimorphisms,  it induces, for any $k\in \mathbb N$, a group homomorphism $j^k_A: H^k_{\mathbb C}(A)\rightarrow H^k_{\mathbb D}(A)$.

\begin{prop}
The homomorphism $j^1_A: H^1_{\mathbb C}(A)\rightarrow H^1_{\mathbb D}(A)$ is an isomorphism.
\end{prop}
\proof
This homomorphism $j^1_A$ is injective: let $X$ be a commutative object in $\mathbb C$ with global support and direction $A$; the image of its class by $j^1_A$ is $0$ when the object $X$ is pointed in $\mathbb D$, which is the case if and only if it is pointed in $\mathbb C$. Let us show $j^1_A$ is surjective: let $Y$ be a commutative object in $\mathbb D$ with global support and direction $A$. Then consider the following pullback of split epimorphisms:
$$\xymatrix@=20pt
{
 Y\times Y \ar@{->>}[r]^{\nu} \ar@<1ex>[d]^{p_{1}} \ar@<-1ex>[d]_{p_{0}} & A \ar[d] \\ 
 Y   \ar@{->>}[r]_{} & 1 \ar@<-1ex>[u]_{e_A}
}
$$
Since $A$ is in $\mathbb C$, then,  by Proposition \ref{n2}, the map $p_0:Y\times Y\rightarrow Y$ is $I$-trivial and these pullbacks are preserved by $I$. They produces the following diagram where any of the upper squares is a pullback:
$$\xymatrix@=20pt
{
 Y\times Y \ar@<2ex>@{->>}[rr]^{\nu} \ar[r]_{} \ar@<1ex>[d]^{p_{1}} \ar@<-1ex>[d]_{p_{0}} & I(Y\times Y) \ar[r]_{I\nu} \ar@<1ex>[d]^{p_{1}} \ar@<-1ex>[d]_{p_{0}} & A \ar[d] \\ 
 Y   \ar@<-2ex>@{->>}[rr]_{} \ar@{->>}[r]^{\eta_Y} \ar@{->>}[d] & IY  \ar@{->>}[r]_{} \ar@{->>}[d]& 1 \ar@<-1ex>[u]_{e_A}\\
 1 \ar@{=}[r]& 1
}
$$
We have $I(Y\times Y)=IY\times IY$. Consequently the two left hand side vertical diagrams are kernel equivalence relations. Since the upper left hand side squares are pullbacks, the Barr-Kock theorem implies that the lower square is a pullback, and that consequently the map $\eta_Y$ is an isomorphism, and $Y$ is in $\mathbb C$.
\endproof

\begin{prop}
The homomorphism $j^2_A: H^2_{\mathbb C}(A)\rightarrow H^2_{\mathbb D}(A)$ is injective.
\end{prop}
\proof
Let $\underline X_1$ be an aspherical groupoid in $\mathbb C$ with direction $A$. Saying that its image by $j^2_A$ is $0$ is saying that there is  a $\underline X_1$-torsor $Y$ in $\mathbb D$, namely a discrete fibration:
$$\xymatrix@=20pt
{
 Y\times Y \ar@{->>}[r]^{\phi_1} \ar@<1ex>[d]^{p_{1}} \ar@<-1ex>[d]_{p_{0}} & X_1 \ar@<1ex>[d]^{d_{1}} \ar@<-1ex>[d]_{d_{0}} \\ 
 Y   \ar@{->>}[r]_{\phi_0}  \ar[u] \ar@{->>}[d]_{} & X_0 \ar[u]\\
 1
}
$$
Again (Proposition \ref{n2}) these pullbacks are preserved by $I$ and produces the following diagram:
$$\xymatrix@=20pt
{
 Y\times Y \ar@<2ex>@{->>}[rr]^{\phi_1} \ar[r]_{} \ar@<1ex>[d]^{p_{1}} \ar@<-1ex>[d]_{p_{0}} & IY\times IY \ar[r]_{I\phi_1} \ar@<1ex>[d]^{p_{1}} \ar@<-1ex>[d]_{p_{0}} & X_1 \ar@<1ex>[d]^{d_{1}} \ar@<-1ex>[d]_{d_{0}} \\ 
 Y   \ar@<-2ex>@{->>}[rr]_{\phi_0} \ar[r]^{\eta_Y}\ar[u] & IY \ar[r]^{I\phi_0}\ar[u] & X_0 \ar[u]
}
$$
Accordingly the left hand side squares are pullbacks, and since $Y$ and $IY$ have global support, the map $\eta_Y$ is an isomorphism for the same reasons as in the previous proposition. Consequently $Y$ is in $\mathbb C$, and the groupoid $\underline X_1$ is in the class of $0$ in $H^2_{\mathbb C}(A)$.
\endproof

\emph{Saying that the homomorphism $j^2_A$ is surjective} is saying that any aspherical groupoid $\underline Y_1$ in $\mathbb D$ with direction $A$ is in the same component class as an aspherical groupoid $\underline T_1$ in $\mathbb C$ with direction $A$, namely that there is a pair of $()_0$-cartesian maps with $\underline Z_1$ aspherical (see \cite{B1}):
$$\underline Y_1 \stackrel{\underline{\theta_1}}{\longleftarrow} \underline Z_1 \stackrel{\underline{\psi_1}}{\longrightarrow} \underline T_1$$
Since $Z_1\twoheadrightarrow Z_0\times Z_0$ is a regular epimorphism, $\underline{\psi_1}$ $()_0$-cartesian and $\underline T_1$ in $\mathbb C$, then, according to Proposition \ref{n2}, the map $Z_1\twoheadrightarrow Z_0\times Z_0$ is $I$-trivial, and the groupoid $\underline Z_1$ is $I$-specific. So, saying that $j^2_A$ is surjective is saying that, for any aspherical groupoid $\underline Y_1$ in $\mathbb D$ with direction $A$, there is a $I$-specific groupoid $\underline Z_1$ in $\mathbb D$ together with a $()_0$-cartesian map $\underline{\theta_1}: \underline Z_1 \rightarrow \underline Y_1$ (which implies that the direction of $\underline Z_1$ is $A$). Accordingly, saying that $j^2_A$ is surjective for any abelian group object $A$ is saying that, for any aspherical groupoid $\underline Y_1$, there is a $I$-specific groupoid $\underline Z_1$ in $\mathbb D$ together with a $()_0$-cartesian map $\underline{\theta_1}: \underline Z_1\rightarrow \underline Y_1$.

\begin{defi}
We shall say that the category $\mathbb D$ is $I$-specific, when, for any aspherical groupoid $\underline Y_1$, there is a $I$-specific groupoid $\underline Z_1$ in $\mathbb D$ together with a $()_0$-cartesian map $\underline{\theta_1}: \underline Z_1 \rightarrow \underline Y_1$.
\end{defi}

With this definition, our previous discussion about the surjectivity of $j^2_A$ becomes:

\begin{prop}\label{n31}
Suppose $\mathbb D$ is an efficiently regular Mal'cev category and $I$ is a reg-epi reflection  satisfying Conditions 1) and 2) of Section \ref{sec2}. Then, the group homomorphism $j^2_A$ is an isomorphism for any abelian group $A$ if and only if the category $\mathbb D$ is $I$-specific.
\end{prop}
Now we get:
\begin{prop}\label{Isp}
Suppose $\mathbb D$ is $I$-specific, then for any aspherical $n$-groupoid $\underline Y_n$ in $\mathbb D$, there is a $I$-specific $n$-groupoid $\underline Z_n$ in $\mathbb D$ together with a $()_{n-1}$-cartesian map $\underline{\theta_n}:\underline Z_n\rightarrow \underline Y_n$.
\end{prop}
\proof
By induction. We suppose the result holds as far the level $n$-$1$. Let $\underline Y_n$ be any aspherical $n$-groupoid and $\theta_{n-1}: \underline Z_{n-1} \rightarrow \underline Y_{n-1}$ make explicit the result at the level $n$-$1$, with $\underline Z_{n-1}$ $I$-specific and $\underline{\theta_{n-1}}$ $()_{n-2}$-cartesian. Then take $\underline{\theta_n}: \underline Z_n \rightarrow \underline Y_n$ the $()_{n-1}$-cartesian map above $\underline{\theta_{n-1}}$. By definition, the $(n$-$1)$-groupoid $\underline Z_{n-1}$ being $I$-specific, so is the $n$-groupoid $\underline Z_n$.
\endproof

\begin{theo}
Suppose $\mathbb D$ is $I$-specific and the assumptions of Proposition \ref{n31} hold. For any abelian group $A$ in $Ab\mathbb C = Ab\mathbb D$, the group homomorphisms $j^n_A: H^n_{\mathbb C}(A)\rightarrow H^n_{\mathbb D}(A)$ are isomorphisms.
\end{theo}
\proof
Let us show that $j^n_A$ is injective. Let $\underline X_n$ be any aspherical $n$-groupoid in $\mathbb C$ with direction $A$ such that its image by $j^n_A$ is $0$. This means that there is a $\underline X_n$-torsor $\underline Y_{n-1}$ in $\mathbb D$, namely a discrete fibration:
$$\xymatrix@=20pt
{
 Y^{\bumpeq}_{n-1} \ar[r]^{\phi_n} \ar@<1ex>[d]^{p_{1}} \ar@<-1ex>[d]_{p_{0}} & X_n \ar@<1ex>[d]^{d_{1}} \ar@<-1ex>[d]_{d_{0}} \\ 
 Y_{n-1}   \ar[r]_{\phi_{n-1}} \ar[u] \ar@{->>}[d]_{(d_0,d_1)} & X_{n-1} \ar[u]\\
 Y^{\bumpeq}_{n-2}
}
$$
Since $\mathbb D$ is $I$-specific, we are able to complete the previous diagram with a $I$-specific $(n-1)$-groupoid $\underline Z_{n-1}$ and a diagram where all the squares are pullbacks:
$$\xymatrix@=20pt
{
 Z^{\bumpeq}_{n-1} \ar[r]^{\theta^{\bumpeq}_{n-1}} \ar@<1ex>[d]^{p_{1}} \ar@<-1ex>[d]_{p_{0}}  &Y^{\bumpeq}_{n-1} \ar[r]^{\phi_n} \ar@<1ex>[d]^{p_{1}} \ar@<-1ex>[d]_{p_{0}} & X_n \ar@<1ex>[d]^{d_{1}} \ar@<-1ex>[d]_{d_{0}} \\ 
 Z_{n-1}   \ar[r]_{\theta_{n-1}} \ar[u] \ar@{->>}[d]_{(d_0,d_1)} & Y_{n-1}   \ar[r]_{\phi_{n-1}} \ar[u] \ar@{->>}[d]_{(d_0,d_1)} & X_{n-1} \ar[u]\\
 Z^{\bumpeq}_{n-2}  \ar[r]_{\theta^{\bumpeq}_{n-2}} & Y^{\bumpeq}_{n-2}
}
$$
Since $\underline X_n$ is in $\mathbb C$ we have the following factorizations where all the squares are pullbacks (again Proposition \ref{n2}):
$$\xymatrix@=20pt
{
 Z^{\bumpeq}_{n-1} \ar@{->>}[r]^{} \ar@<1ex>[d]^{p_{1}} \ar@<-1ex>[d]_{p_{0}}  &I(Z^{\bumpeq}_{n-1}) \ar[rr]^{I(\phi_n.\theta^{\bumpeq}_{n-1})} \ar@<1ex>[d]^{p_{1}} \ar@<-1ex>[d]_{p_{0}} & & X_n \ar@<1ex>[d]^{d_{1}} \ar@<-1ex>[d]_{d_{0}} \\ 
 Z_{n-1}   \ar@{->>}[r]_{\eta_{n-1}} \ar[u] \ar@{->>}[d]_{(d_0,d_1)} & I(Z_{n-1})   \ar[rr]_{I(\phi_{n-1}.\theta_{n-1})} \ar[u] \ar@{->>}[d]_{(d_0,d_1)} & & X_{n-1} \ar[u]\\
 Z^{\bumpeq}_{n-2}  \ar@{->>}[r]_{} & I(Z^{\bumpeq}_{n-2})
}
$$
Now, since $\underline Z_{n-1}$ is $I$-specific, we have certainly $I(Z^{\bumpeq}_{n-1})=I(Z_{n-1})^{\bumpeq}$ and $I(Z^{\bumpeq}_{n-2})=I(Z_{n-2})^{\bumpeq}$, and thus a $\underline X_n$-torsor in $\mathbb C$:
$$\xymatrix@=20pt
{
 I(Z_{n-1})^{\bumpeq} \ar[rrr]^{I(\phi_n.\theta^{\bumpeq}_{n-1})} \ar@<1ex>[d]^{p_{1}} \ar@<-1ex>[d]_{p_{0}} &&& X_n \ar@<1ex>[d]^{d_{1}} \ar@<-1ex>[d]_{d_{0}} \\ 
 I(Z_{n-1})   \ar[rrr]_{I(\phi_{n-1}.\theta_{n-1})} \ar[u] \ar@{->>}[d]_{(d_0,d_1)} &&& X_{n-1} \ar[u]\\
 I(Z_{n-2})^{\bumpeq}
}
$$
Let us show now that $j^n_A$ is surjective. Let $\underline Y_n$ be any aspherical $n$-groupoid in $\mathbb D$ with direction $A$. According to Proposition \ref{Isp}, there is a $I$-specific $n$-groupoid $\underline Z_n$ in $\mathbb D$ together with a $()_{n-1}$-cartesian map $\underline{\theta_n}:\underline Z_n\rightarrow \underline Y_n$. Then the following pair of maps shows that the image by $j^n_A$ of $I\underline Z_n$ (which is in $\mathbb C$) is in the same component class as $\underline Y_n$:
$$I\underline Z_n \stackrel{\underline{\eta_n}}{\longleftarrow} \underline Z_n \stackrel{\underline{\theta_n}}{\longrightarrow} \underline Y_n$$
since $\underline{\eta_n}$ is $()_{n-1}$-cartesian by Proposition \ref{ff}.
\endproof

\section{Examples of $I$-specific categories}

\subsection{The protomodular context}

Given any finitely complete category $\mathbb E$, we denote by $Pt\mathbb E$ the category whose objects are the split epimorphisms in $\mathbb E$ and whose morphisms are the commutative squares between them:
$$\xymatrix@=20pt
{
 X' \ar[r]^{x} \ar@<-1ex>[d]_{f'} & X \ar@<-1ex>[d]_{f}  \\ 
 Y' \ar[r]_{y} \ar[u]_{s'} & Y \ar[u]_{s}  
}
$$
The codomain functor: $Pt\mathbb E\rightarrow \mathbb E$ is a fibration which is called \emph{the fibration of points} and whose cartesian maps are those previous squares which are pullbacks. The fibre $Pt_Y\mathbb E$ above $Y$ has the split epimorphisms above $Y$ as objects and the commutative triangles between them as morphisms. Recall that the category $\mathbb E$ is \emph{protomodular} \cite{B5} when the change of base functors with respect to this fibration are conservative (i.e. reflect the isomorphisms) and that any protomodular category is necessarily a Mal'cev category \cite{BB}. The category $Gp$ of group, the category $Rg$ of non unitary commutative rings, the category $\mathbb K$-$Lie$ of $\mathbb K$-Lie algebras and any additive category are pointed protomodular. It is the case also of the category $GpTop$ of topological groups, and more generally of any category $Top^{\mathbb T}$ of topological protomodular algebras \cite{BC}. 

Given any finitely complete category $\mathbb E$, any fibre of the fibration $()_0: Grd\mathbb E\rightarrow \mathbb E$ is protomodular. When moreover $\mathbb E$ is protomodular, any slice category $\mathbb E/Y$ is protomodular, while any fibre $Pt_Y\mathbb E$ is pointed protomodular. It appears that the context of pointed protomodular categories particularly fits with the treatment of exact sequences \cite{B5}; in particular the (split) short five lemma holds in it \cite{BB}.

\subsection{Peri-abelian categories}

Consider first a finitely cocomplete regular pointed protomodular category $\mathbb D$ and $\mathbb C=Ab\mathbb D$ the subcategory of the abelian (group) objects (actually, $\mathbb D$ being pointed, we have $Ab\mathbb D=M\mathbb D$ since any commutative object, being pointed, is abelian). Then the inclusion $Ab\mathbb D\rightarrowtail \mathbb D$ has a reg-epi reflection which will be denoted by $A$. This reflection obviously satisfies Condition 1) of Section \ref{sec2}. It satisfies Condition 2) since, the category $\mathbb D$ being pointed, it preserves any product by Proposition \ref{pb}. More generally, for any object $Y$ in $\mathbb D$, any inclusion $Ab(Pt_Y\mathbb D)\rightarrowtail Pt_Y\mathbb D$ produces a reg-epi reflection denoted by $A_Y$, since  the category of $Pt_Y\mathbb D$ of split epimorphisms in $\mathbb D$ above $Y$ is still finitely cocomplete, regular, and pointed protomodular, and this reflection $A_Y$ satisfies Conditions 1) and 2). In the case $\mathbb D=Gp$ the category of groups, the functor $A_Y$ is classically described by the following diagram where $\bar X=X/[\ker f ,\ker f]$, since $[\ker f ,\ker f]$ is a normal subgroup of $X$: 
$$\xymatrix@=15pt
{
 X  \ar[rd]^f \ar@{->>}[rr]^{\eta_f}&& \bar X \ar[ld]^{\bar f}\\
 & Y \ar@<1ex>[lu]^{s} \ar@<1ex>[ru]
}
$$
It is clear that $ker{\bar f}=\ker f/[\ker f ,\ker f]=A(\ker f)$, so that the change of base with respect to the fibration of points along the initial map $\alpha_Y: 1\rightarrow Y$ preserves the associated abelian object. Consequently this is the case for the change of base along any map $Y'\rightarrow Y$. This property is far from being satisfied by any finitely cocomplete regular pointed protomodular or even semi-abelian category (i.e. non only regular, but also exact, see \cite{JMT}) as shown by the next lemma below.

\begin{defi}\label{peri}
We shall say that a finitely cocomplete, regular, pointed protomodular category $\mathbb D$ is peri-abelian when the change of base functor along any map $Y'\rightarrow Y$ with respect to the fibration of points preserves the associated abelian object.
\end{defi}

It is clear that this definition still holds in the finitely cocomplete, regular, non-pointed Mal'cev context, but we shall only focus here on the context precised in the definition. So, here, ''peri-abelian" will imply ''finitely cocomplete, regular, pointed protomodular". If $(AbPt)\mathbb D$ denotes the subcategory of the abelian objects in the fibres of the fibration of points, it is equivalent to saying that the reg-epi reflection $A_{()}$ is cartesian, i.e. it preserves the cartesian maps:
$$\xymatrix@=10pt
{
{(AbPt)\mathbb D\;\;} \ar[rdd] \ar@{>->}[rr]^{}&& Pt\mathbb D \ar[ldd] \ar@(u,u)[ll]_{A_{()}}\\
\\
 & \mathbb D 
}
$$
For the same reasons as for the category $Gp$, it is clear that the categories $Rg$ of non unitary commutative rings and $\mathbb K$-$Lie$ of $\mathbb K$-Lie algebras are peri-abelian. When $X$ is a topological group, endowing the normal subgroups $\ker f$ and $[\ker f ,\ker f]$ with the induced topology allows us to use the same argument as above to show that the category $GpTop$ of topological groups is peri-abelian.

\begin{lemma}
Suppose $\mathbb D$ is peri-abelian. Then the change of base functors with respect to the fibration of points along any map {do \rm reflect} the abelian objects. 
\end{lemma}
\proof
It is sufficient to prove it for the initial map $\alpha_Y: 1\rightarrowtail Y$. So, suppose $\ker f=A$ is abelian, where $:X\rightarrow Y$ is a split epimorphism. Since $\mathbb D$ is peri-abelian, we have also $\ker{\bar f}=A$; more precisely, the factorization determined by $\eta_f$ is an isomorphism. So, by the short five lemma, the map $\eta_f$ is itself an isomorphism, and $f$ is in $(AbPt)\mathbb D$.
\endproof

Accordingly the semi-abelian category $DiGp$ of digroups (a digroup is a set endowed with two group structures whose only coherence condition is to have same unit element) is not peri-abelian since its change of base functors with respect to the fibration of points do not reflect the abelian objects.\\

We are going now to extend the preservation of abelian objects at the level of the fibration of points to the preservation of commutative objects at the level of the slice categories.

Suppose $\mathbb D$ is not only regular, but efficiently regular. Consider the slice category $\mathbb D/Y$. Then the construction of the associated commutative object in $\mathbb D/Y$ can be derived from the  construction of the associated abelian object in the fibres $Pt_U\mathbb D$. For that, given any map $f:X\rightarrow Y$, consider the following diagram, where the two lower vertical maps $\bar p_0$ are the associated abelian objects of the pointed objects $p_0:R[f]^2\rightarrow R[f]$ and $p_0:R[f]\rightarrow X$ respectively :
$$\xymatrix@=20pt
{
  R^2[f]  \ar@<1ex>[r]^{p_{1}} \ar@<-1ex>[r]_{p_{2}} \ar@{->>}[d] \ar@<-3ex>[dd]_{p_0} & R[f] \ar@{->>}[r]^{d_1} \ar@{->>}[d]  & X \ar@{->>}[d]^{\eta} \ar@(r,r)[dd]^f\\
  \overline{R[f]^2}  \ar@<1ex>[r]^{\overline{p_{1}}} \ar@<-1ex>[r]_{\overline{p_{2}}} \ar[d]_{\bar{p_0}} & \overline{R[f]} \ar@{.>>}[r]^{\bar q} \ar[d]_{\bar{p_0}} & \bar X \ar[d]^{\bar f}\\
  R[f]  \ar@<1ex>[r]^{p_{0}} \ar@<-1ex>[r]_{p_{1}} & X  \ar[r]_{f}  & Y 
}
$$ 
Since the left hand side vertical rectangles are pullbacks of split epimorphisms, and the two upper vertical arrows are regular epimorphisms, then, as shown in \cite{BR2}, in the regular Mal'cev category $\mathbb D$, the two left hand side levels of squares are pullbacks, and produce two discrete fibrations. In particular, since the lower one is a discrete fibration above an equivalence relation, the middle horizontal diagram is an equivalence relation. Since $\mathbb D$ is efficiently regular, it is an effective equivalence relation. Let $\bar q$ be its coequalizer and $\bar f$, $\eta$  be the induced factorizations. Now since all the squares on the left hand side are pullbacks, and the maps $d_1$ and $\bar q$ are regular epimorphisms, the two squares on the right hand side are pullbacks. Moreover, the map $\bar f$ has an abelian kernel equivalence relation as a quotient of the map $\bar p_0$ which has an abelian kernel relation. It is straighforward to show that $\bar f$ is the associated abelian object of $f$ from the fact that $\bar p_0$ is the associated abelian object to $p_0:R[f]\rightarrow X$.

\begin{prop}
When $\mathbb D$ is efficiently regular and peri-abelian, then, given any map $h: Y'\rightarrow Y$ in $\mathbb D$, the change of base functor $h^*:\mathbb D/Y\rightarrow \mathbb D/Y'$ preserves the associated commutative object. 
\end{prop}
\proof
Consider the following diagram were the lower (and thus any) square is a pullback:
$$\xymatrix@=20pt
{
 R[f'] \ar[r]^{R(g)} \ar@<1ex>[d]^{p_{1}} \ar@<-1ex>[d]_{p_{0}} & R[f]\ar@<1ex>[d]^{p_{1}} \ar@<-1ex>[d]_{p_{0}}\\
 X'\ar[r]^g \ar[d]_{f'}& X \ar[d]^{f}\\
 Y'\ar[r]_h & Y 
}
$$
Then, with the notation of the previous construction, the two following vertical rectangles are the same:
$$\xymatrix@=20pt
{
 \overline{R[f']} \ar[r]^{\overline{R(g)}} \ar[d]_{\overline{p_{0}}}  & \overline{R[f]}\ar[d]^{\overline{p_{0}}} &&&   \overline{R[f']} \ar[r]^{\overline{R(g)}}  \ar@{->>}[d]_{q_{X'}}  & \overline{R[f]} \ar@{->>}[d]^{q_{X}}\\
 X'\ar[r]^g \ar[d]_{f'}& X \ar[d]^{f} &&& \bar X' \ar[r]^{\bar g} \ar[d]_{\bar{f'}} & \bar X \ar[d]^{\bar f}\\
 Y'\ar[r]_h & Y  &&& Y'\ar[r]_{h} & Y
}
$$
The lower left hand side square is a pullback by assumption, while the upper left hand side square is a pullback since $\mathbb D$ is peri-abelian. Accordingly the right hand side rectangle is a pullback. If the upper right hand side square is shown to be a pullback, the lower right hand side square will be a pullback since the two upper vertical arrows are regular epimorphism. Now consider the following horizontal rectangle made of two pullbacks:
$$\xymatrix@=20pt
{
  \overline{R[f']} \ar[r]^{\bar p_0} \ar@{->>}[d]_{q_{X'}} & X'\ar[r]^g \ar[d]_{f'}& X \ar[d]^{f}\\
  \bar X' \ar[r]_{\bar{f'}} & Y'\ar[r]_h & Y
}
$$
it is the same as the following one:
$$\xymatrix@=20pt
{
  \overline{R[f']} \ar[r]^{\overline{R(g)}} \ar@{->>}[d]_{q_{X'}} & \overline{R[f]}\ar@{->>}[d]^{q_{X}}\ar[r]^{\bar p_0}& X \ar[d]^{f}\\
  \bar X' \ar[r]_{\bar g} &\bar X \ar[r]_{\bar f} & Y
}
$$
Now, the whole rectangle and the right hand side square are pullbacks, accordingly the left hand side square is a pullback.
\endproof

\subsection{The projective assumption}\label{sec}

We shall suppose, unless otherwise stated, that the pointed category $\mathbb D$ is finitely cocomplete, efficiently regular, protomodular and peri-abelian. This holds in particular when $\mathbb D$ is semi-abelian and peri-abelian, as are the categories $Gp$, $Rg$ and $\mathbb K$-$Lie$. This holds also when $D$ is the category $GpTop$ of topological groups. Under our conditions, it is clear that the reg-epi reflections $A:\mathbb D\rightarrow Ab\mathbb D$, $A_Y:Pt_Y\mathbb D\rightarrow AbPt_Y\mathbb D$ and $M_Y:\mathbb D/Y\rightarrow M(\mathbb D/Y)$ satisfy  Condition 1) of Section \ref{sec2}. We already noticed that the two first ones satisfy Condition 2). It is shown in \cite{BR2} that, in presence of the efficiently regular assumption, the third one satisfies also Condition 2).

As usual, we have the following:
\begin{defi}
Given any homological (i.e. pointed regular and protomodular) category $\mathbb E$, we shall say that the object $Y$ has projective dimension $1$ when there is an exact sequence with $H$ and $K$ projective objects in $\mathbb E$ (with respect to the regular epimorphisms):
$$\xymatrix@=20pt
{
 1 \ar[r] & {K\;} \ar@{>->}[r]^ k & H \ar@{->>}[r]^ h & Y \ar[r] & 1 
}
$$
\end{defi}

We are now going to prove that, \emph{when the object $Y$ in our category $\mathbb D$ has projective dimension $1$, the slice category $\mathbb D/Y$ is $M_Y$-specific}, where $M_Y$ is the reg-epi reflection associated with the inclusion of the commutative objects $M(\mathbb D/Y) \rightarrowtail \mathbb D/Y$. For that, we shall need some recall about the construction of the normalization of a groupoid.

\subsection{The normalization functor}

The category $\mathbb D$ being pointed, we can associate with any internal groupoid $\underline U_1$ its Moore normalization, namely the 
$1$-chain complex (= map) given by the left hand side vertical arrow in the following pullback:
$$
\xymatrix@=20pt
{
  K[d_0] \ar[r]^{\kappa} \ar[d]_{\mu_{\underline U_1}} & U_1  \ar[d]^{(d_0,d_1)} \\
  U_0 \ar[r]_{(0,1_{U_0})} & U_0\times U_0
}
$$
Clearly this construction extends to a functor: $\mu: Grd\mathbb D\rightarrow Ch_{1}\mathbb D$ which maps $()_0$-cartesian functors onto pullback squares. For us, the main fact will be that, since $\mathbb D$ is efficiently regular and protomodular, the restriction of this functor $\mu$ to the subcategory of aspherical (=connected since $\mathbb D$ is pointed) groupoids determines an equivalence of categories onto the category of central extensions in $\mathbb D$ (see Proposition 17 in \cite{B5}). So, let us consider any aspherical groupoid in $\mathbb D/Y$:
$$
\xymatrix@=30pt
{ 
  H_1 \ar@<2ex>[r]^{d_1} \ar@<-2ex>[r]_{d_0} & H_0 \ar[l]_{s_0} \ar@{->>}[r]^{h} & Y
 }
$$ 
\begin{lemma}
Suppose $\mathbb D$ is a pointed regular and protomodular category. If the kernel $K$ of the regular epimorphism $h$ is projective in $\mathbb D$, then the aspherical groupoid $\alpha_Y^*(\underline H_1)$ in $\mathbb D$ is pointed. If $\mathbb D$ is finitely cocomplete, the aspherical groupoid $\alpha_Y^*(\underline H_1)$ is $A$-specific. If moreover $\mathbb D$ is efficiently regular and peri-abelian, any aspherical groupoid $\underline H_1$ in $\mathbb D/Y$ (above $h$) is $M_Y$-specific.
\end{lemma}
\proof
Consider the following diagram where any square is a pullback:
$$
\xymatrix@=20pt
{ 
  {N\underline K_1\;} \ar@<2ex>[rr]^{\kappa}  \ar@{>->}[r]_{\bar{\kappa}} \ar@{->>}[d]_{\mu_{\underline K_1}} & {K_1\;} \ar@{>->}[r]^{k_1} \ar@{->>}[d]^{(d_0,d_1)}& H_1 \ar@{->>}[d]^{(d_0,d_1)} \\
  K \ar[d]  \ar@{>->}[r]^{(0,1)} & {K\times K \;} \ar@<1ex>[d]^{p_{1}} \ar@<-1ex>[d]_{p_{0}}\ar@{>->}[r]^{R(k)} & R[h] \ar@<1ex>[d]^{p_{1}} \ar@<-1ex>[d]_{p_{0}}\\
  1  \ar[r] & {K\;}\ar[d] \ar@{>->}[r]^{k} & H_0 \ar@{->>}[d]^h\\ 
  & 1 \ar[r]_{\alpha_Y} & Y
 }
$$
The groupoid $\alpha_Y^*(\underline H_1)$ is given by the middle vertical diagram, and denoted by $\underline K_1$. Its normalization $\mu_{\underline K_1}$ is a regular epimorphism, as a pullback of a regular epimorphism. The object $K$ being projective, this regular epimorphism $\mu_{\underline K_1}$ is split. Since the groupoid $\underline K_1$ is aspherical and the restriction of the normalization functor $\mu$ to the aspherical groupoids produces an equivalence of categories with the central extensions in $\mathbb D$, the splitting of $\mu_{\underline K_1}$ determines a splitting of the core $(d_0,d_1): K_1\twoheadrightarrow K\times K$ which makes the groupoid $\underline K_1$ a pointed groupoid. According to Example \ref{pointed}, the fact that $\underline K_1$ is pointed makes it $A$-specific where $A: \mathbb D\rightarrow Ab\mathbb D$ give the associated abelian object. This implies that $\underline{\eta}_1: \underline K_1\rightarrow A(\underline K_1)$ is  $()_0$-cartesian in $\mathbb D$. 

When, moreover, $\mathbb D$ is efficiently regular and peri-abelian, the functor $\alpha_Y^*: \mathbb D/Y\rightarrow \mathbb D$ preserves the associated commutative objects, and we have $\underline{\eta}_1=\alpha_Y^*(\underline{\eta}_1^Y)$, with $\underline{\eta}_1^Y:\underline H_1 \rightarrow M_Y(\underline H_1)$  in $\mathbb D/Y$. It remains to show this last map is $()_0$-cartesian in $\mathbb D/Y$. Since pulling back along $\alpha_Y: 1\rightarrowtail Y$ determines a functor $\alpha_Y^*:\mathbb D/Y\rightarrow \mathbb D$ which reflects isomorphisms, when it is restricted to regular epimorphisms of $\mathbb D/Y$ (this is the short five lemma), the pullback making the functor $\underline{\eta}_1$ a $()_0$-cartesian morphism in $\mathbb D$ is reflected in $\mathbb D/Y$. Accordingly the functor $\underline{\eta}_1^Y$ is $()_0$-cartesian in $\mathbb D/Y$ and the groupoid $\underline H_1$ is $M_Y$-specific.
\endproof

Whence the following:
\begin{theo}
Under the conditions on $\mathbb D$ assumed at the beginning of Section \ref{sec}, when the object $Y$ has projective dimension $1$, the category $\mathbb D/Y$ is $M_Y$-specific.
\end{theo}
\proof
Keeping the notations of the previous definition, let $Y$ have projective dimension $1$. Let us consider any internal aspherical groupoid in $\mathbb D/Y$:
$$
\xymatrix@=30pt
{ 
  Y_1 \ar@<2ex>[r]^{d_1} \ar@<-2ex>[r]_{d_0} & Y_0 \ar[l]_{s_0} \ar@{->>}[r]^{q} & Y
 }
$$
Since the domain $H$ of $h:H\twoheadrightarrow Y$ is projective, there is a factorization $\theta :H\rightarrow Y_0$ such that $q.\theta=h$. Let $\underline{\theta_1}: \underline H_1 \rightarrow \underline Y_1$ be the $()_0$-cartesian map above $\theta$. The map $h:H\twoheadrightarrow Y$ being a regular epimorphism, it is clear that the groupoid $\underline H_1$ is aspherical in the category $\mathbb D/Y$. Since, moreoever, the kernel $K$ of $h$ is projective, according to the previous lemma, the  groupoid $\underline H_1$  produces, inside the category $\mathbb D/Y$, the $M_Y$-specific groupoid associated with $\underline Y_1$ we are looking for. 
\endproof
As a straightforward consequence, we get:
\begin{theo}
Under the conditions on $\mathbb D$ assumed at the beginning of Section \ref{sec}, when the object $Y$ has projective dimension $1$, then, for any abelian group $A$ in the slice category $\mathbb D/Y$ and any integer $n$, the cohomology groups $H^n_{M(\mathbb D/Y)}(A)$ and $H^n_{\mathbb D/Y}(A)$ are isomorphic.
\end{theo}

It is well known that any subgroup of a free group is free. Accordingly, given any group $Y$, the following exact sequence makes explicit the projective dimension $1$ of the group $Y$, where $FU(Y)$ is the free group above the underlying set $U(Y)$:
$$\xymatrix@=20pt
{
 1 \ar[r] & {K\;} \ar@{>->}[r]^ k & FU(Y) \ar@{->>}[r] & Y \ar[r] & 1 
}
$$
So that, for any group $Y$, the previous theorem holds. Our groups $H^n_{Gp/Y}(A)$ are nothing but the ones given by Eilenberg-Mac Lane cohomology, since the category of crossed $n$-fold exensions with kernel $A$ by means of which these cohomology groups are realized in \cite{Ho} and \cite {Hu}, coincides with the category of aspherical $n$-groupoids with direction $A$ in the slice category $Gp/Y$ (see \cite{B1}). The translation of our theorem means that any class of crossed $n$-fold exensions, for a given $Y$-module structure on the abelian group $A$:
$$ 1\rightarrow A\rightarrow Y_n  \;\;...\;\; Y_1\rightarrow Y_0\stackrel{q}{\rightarrow} Y \rightarrow 1$$ 
has an equivalent representation where $Z_1$ is abelian:  
$$ 1\rightarrow A\rightarrow Z_n  \;\;...\;\; Z_1\rightarrow Z_0\stackrel{\bar q}{\rightarrow} Y \rightarrow 1$$
This was first made explicit in Proposition 2.7 in the Holt's article about the Eilenberg Mac-Lane cohomology groups \cite{Ho}, while the Huebschmann's method was quite different \cite{Hu}. In a way, our isomorphisms show that, via the ''peri-abelian" connection and the projective dimension, the categories $Gp$ and $Ab$ are so close to each other that the Eilenberg-Mac Lane cohomology is unable to discriminate them. 

On the other hand the so-called Shirshov-Witt theorem asserts that any Lie subalgebra of a free Lie algebra is free. Accordingly the previous theorem also holds in the category $\mathbb K$-$Lie$, for any $\mathbb K$-Lie algebra $Y$.\\

Since the category $GpTop$ of topological groups satisfies all the conditions given at the beginning of Section \ref{sec}, a natural question would be wether it is possible to adapt, or not, the last step of our method (dealing with the projective dimension) to obtain similar cohomology isomorphisms in this topological setting.

\bigskip\bigskip\noindent{\it Universit\'e du Littoral, Calais, France\\ bourn@lmpa.univ-littoral.fr}

\end{document}